\documentclass [12pt,a4paper]{amsart}

\usepackage{amssymb,amsxtra,amsfonts}
\usepackage{lscape}

\openup\jot 
\newcommand\beq{\begin{equation}}
\newcommand\eeq{\end{equation}}

\newcommand{\IP}{\mathbb{P}}                                     
                           
\newcommand{\IR}{\mathbb{R}}                           
\newcommand{\IC}{\mathbb{C}}

\newcommand{\cC}{\mathcal{C}}

\newcommand{\cH}{\mathcal{H}}

\newcommand{\PVI}{$\text{P}_{\text{VI}}$}   %

\newcommand{\lsl}{\mathfrak{sl}} 

\newcommand{\bfa}{{\bf a}}
\newcommand{\bfb}{{\bf b}}
\newcommand{\bfc}{{\bf c}}
\newcommand{\bfd}{{\bf d}}
\newcommand{\bfe}{{\bf e}}
\newcommand{\bff}{{\bf f}}

\newcommand{\pf}{\begin{bpf}}

\newcommand{\pfms}{\begin{bpfms}}
\newcommand{\epf}{\end{bpf}\hfill$\square$\\}           %
\newcommand{\epfms}{\end{bpfms}\hfill$\square$\\}               %

\newcommand{\idea}{\begin{bidea}}

\newcommand{\eidea}{\end{bidea}\hfill$\square$\\}           %

\newcommand{\sk}{\begin{bsk}}    %

\newcommand{\esk}{\end{bsk}\hfill$\square$\\}           %
\newcommand{\sketch}{\begin{bsketch}}%

\newcommand{\esketch}{\end{bsketch}\hfill$\square$\\}

\newcommand{\wh}{\widehat}
\newcommand{\al}{\alpha}

\newcommand{\be}{\beta}
\newcommand{\ga}{\gamma}
\newcommand{\de}{\delta}

\newcommand{\la}{\lambda}

\newcommand{\Sym}{\text{\rm Sym}}

\newcommand{\sym}{\text{\rm Sym}}
\newcommand{\tr}{\text{\rm Tr}}

\newcommand{\Hom}{\text{\rm Hom}}

\newcommand{\SL}{\text{\rm SL}}

\newcommand{\GL}{\text{\rm GL}}

\newcommand{\SU}{\text{\rm SU}}

\theoremstyle{plain}
\newtheorem {hypo}{\bf\hspace{-\parindent}Hypothesis}

\newtheorem {thm}{Theorem}
\newtheorem {prop}[hypo]{Proposition}%

\newtheorem {lem}[hypo]{Lemma}%

\theoremstyle{definition}

\theoremstyle{remark}
\newtheorem {rmk}[hypo]{Remark}%

\begin{document}

%\date{\today}

\title[Regge and Okamoto symmetries]{Regge and Okamoto symmetries}
\author{Philip P. Boalch}
\address{\'Ecole Normale Sup\'erieure\\
45 rue d'Ulm\\
75005 Paris\\
France} 
\email{boalch@dma.ens.fr}
\urladdr{www.dma.ens.fr/$\sim$boalch}

\begin{abstract}
We will relate the surprising
Regge symmetry of the Racah-Wigner $6j$ symbols
to the surprising Okamoto symmetry of the Painlev\'e VI differential equation.
This then presents the opportunity to give a conceptual derivation of the Regge
symmetry, as the representation theoretic analogue of the derivation in
\cite{k2p, srops} of the Okamoto symmetry. 
\end{abstract}

\maketitle

\renewcommand{\baselinestretch}{1.02}            %
\normalsize

\begin{section}{Introduction}

The $6j$-symbols (or Racah coefficients) are real numbers
 associated to the choice of six irreducible representations 
$V_a,\ldots,V_f$ of $\SU(2)$.
They were first published in work of Racah \cite{Racah42} in
1942, and arise in the addition of the three angular momenta,
which classically can be viewed as adding three vectors in $\IR^3$.
Apparently (\cite{PzRg}) ``there is hardly any branch of physics involving
angular momenta where the use of Racah-coefficients is not needed in order to
carry out the simplest computation''.
(See the two volumes \cite{BiedLou1,BiedLou2} for many more details or the
introduction to the tables \cite{Rotetal} for a concise summary.)
Wigner \cite{Wigner40}
used a
slightly different normalisation so that they have 
tetrahedral symmetry, and wrote them in
the form:
\beq\label{eq: 6jsymbol}
\left\{\begin{matrix} a&b&e\\c&d&f\end{matrix} \right\}.
\eeq
Here $a,b,c$ should be thought of as the lengths of three vectors 
$\bfa,\,\bfb,\,\bfc$ in $\IR^3$ so the four points 
${\bf 0},\,\bfa,\,\bfa+\bfb,\,\bfa+\bfb+\bfc$ are 
the vertices of a (skew) tetrahedron.
Then $d,e,f$ should be the lengths of the other three edges of this
tetrahedron, i.e. the lengths of $\bfa+\bfb+\bfc, \,\bfa+\bfb,\, \bfb+\bfc$
respectively. (Thus each column of \eqref{eq: 6jsymbol} contains the 
lengths of two opposite edges, and the top row $abe$ is a face.)
Then the $6j$-symbol is invariant under the possible relabellings
of this tetrahedron (preserving the relations so one gets $24=\#\sym_4$
possibilities).

Racah established an explicit formula for the $6j$-symbols as a sum,
which has since been equated with the value at $1$ of certain $_4F_3$
hypergeometric functions (see e.g. \cite{Wilson-orthog}).

Using this explicit Racah formula, in 1959 Regge \cite{Reg-6j}
showed the $6j$-symbols also have the 
following further symmetry, which is more mysterious:
$$\left\{\begin{matrix} a&b&e\\c&d&f\end{matrix} \right\}
= \left\{\begin{matrix} p-a&p-b&e\\p-c&p-d&f\end{matrix} \right\}
$$
where $p=(a+b+c+d)/2$.
(Combined with the tetrahedral symmetries this generates a symmetry group
isomorphic to $\sym_3\times\sym_4$.)

For example classically, recalling that a tetrahedron
is determined up to isometry by its edge lengths, one may show 
(cf. Ponzano--Regge \cite{PzRg} and Roberts \cite{Roberts-6j1}) that
this Regge action on the set of six edge lengths  
defines a non-trivial
automorphism of the set of Euclidean tetrahedra, taking a generic 
tetrahedron to
a {\em non-congruent} tetrahedron.

Earlier Regge \cite{Reg-3j} 
found similar extra symmetries of the
Clebsch-Gordan $3j$-symbols. The $3j$-symbols 
are in a sense less canonical, but note
that Ponzano--Regge \cite{PzRg} p.7 explain how to obtain the $3j$-symbols as
an asymptotic limit of $6j$-symbols and that in this limit the $6j$ Regge 
symmetry becomes the $3j$ Regge symmetry. They then wrote:
\begin{quote}
{\em The geometrical and physical content of these [Regge] symmetries is
  still to be understood and they remain a puzzling feature of the theory of
  angular momenta. Therefore it is a pleasant result to be able to reduce the
  problem of their interpretation to the Racah coefficient only.} 
\end{quote}

The basic aim of this article is to give a conceptual explanation of the
Regge $6j$ symmetry.
The key idea is to relate the above $6j$-symbols (for the group $\SU(2)$) to
certain three-dimensional $6j$-symbols (i.e. for the group $\SU(3)$). 
The Regge transformation then arises from the natural duality between two
dual irreducible representations of $\SU(3)$. 

The layout of the remainder of this article is as follows.
First we will give the definition of the $6j$-symbols, then we will
relate the Regge symmetry to a symmetry of a
completely different object, this time a non-linear differential equation, the
Okamoto symmetry of the Painlev\'e VI equation. 

Then we will ``quantise'' (that is, give the representation theoretical
analogue of) the derivation of the Okamoto symmetry given in 
\cite{k2p, srops}, and so give a conceptual derivation of the Regge symmetry
(i.e. without using the Racah formula).

\end{section}

\begin{section}{Background}

The $6j$-symbols are real numbers associated to the choice of $6$ irreducible
representations (irreps) of 
$G:=\SU(2)$.
We will label irreps $V_a$ by positive integers $a\ge 0$, so that 
$V_a= \Sym^a(V)$ is the spin $a/2$ representation of dimension $a+1$, 
where $V$ is the two-dimensional Hermitian vector
space defining $G$.
Given $3$ such irreps, with labels $a,b,c$ say, one can form the 
tensor product
$$V_{abc}:=V_a\otimes V_b\otimes V_c$$
which again will be a representation of $G$ and will decompose as a direct sum
of irreps. Thus given a fourth representation $V_d$ one may consider the
multiplicity space
$$M_{abcd}:=\Hom_G(V_d,V_{abc})$$
of $G$-equivariant maps from $V_d$ into the 3-fold tensor product.
Thus $M_{abcd}$ is a Hermitian vector space with dimension 
equal to the multiplicity of 
$V_d$ in $V_{abc}$.

There are two (almost) canonical unitary 
bases in $M_{abcd}$ (`coupling bases') and the $6j$-symbols
arise as matrix entries of the change of basis matrix
between these two bases.
The coupling 
bases arise by decomposing $V_{abc}$ in the two possible orders:
on one hand we may first decompose $V_{ab}:=V_a\otimes V_b$:
$$V_{ab}\cong \bigoplus_e V_e\otimes \Hom_G(V_e,V_{ab})$$
which entails the following direct sum decomposition of $M_{abcd}$:
$$M_{abcd}=\Hom_G(V_d,V_{ab}\otimes V_c)\cong
\bigoplus_e \Hom_G(V_d,V_{ec})\otimes \Hom_G(V_e,V_{ab}).
$$
The key point is that each of the terms in this direct sum is either zero or
one-dimensional (since, for $\SU(2)$, 
any irrep appears at most once in the tensor product of two
irreps).
Thus choosing a real vector of length one in each one dimensional 
term yields the $1$-$2$ coupling  
basis $\{v_{e}\}$  of $M_{abcd}$ as $e$ varies, 
unique up to the sign of each basis vector.
(We set $v_e=0$ if the space $\Hom_G(V_d,V_{ec})\otimes
\Hom_G(V_e,V_{ab})$ is zero.)
Similarly decomposing the 3-fold product in the other order, i.e. first writing
$$V_{bc}\cong \bigoplus_f V_f\otimes \Hom_G(V_f,V_{bc})$$
yields a different basis $\{w_f\}$, the $2$-$3$ coupling basis, 
adapted to the decomposition
$$M_{abcd}=
\bigoplus_f \Hom_G(V_d,V_{af})\otimes \Hom_G(V_f,V_{bc}).
$$
Thus given six irreps with labels $a,b,c,d,e,f$, and a standard 
sign-convention, 
one will get two vectors 
$v_e,w_f$ in $M_{abcd}$ and thus a number
$$U(a,b,c,d,e,f)= \langle v_e,w_f \rangle$$
by pairing them using the Hermitian form.
(As $e$ and $f$ vary these will be the matrix entries of the unitary change of
basis matrix alluded to above---in fact by reality it is real orthogonal.)

The $6j$-symbols were defined by Wigner in 
terms of $U$ by a minor normalisation:
$$\left\{\begin{matrix} a&b&e\\c&d&f\end{matrix} \right\}
= (-1)^{p}\frac{U(a,b,c,d,e,f)}{\sqrt{(e+1)(f+1)}},\quad
p:=(a+b+c+d)/2.$$
This normalisation is such that the $6j$-symbols admit tetrahedral symmetry,
where the coefficients label the six edges of a tetrahedron (containing for
example the
quadrilateral $abcd$ and faces $abe$ and $bcf$).
Note that if $U$ is non-zero then $p$ will be an integer.

As mentioned in the introduction using the explicit Racah formula for the
$6j$-symbols, 
Regge \cite{Reg-6j}
showed the $6j$-symbols also have the 
following symmetry:
$$\left\{\begin{matrix} a&b&e\\c&d&f\end{matrix} \right\}
= \left\{\begin{matrix} p-a&p-b&e\\p-c&p-d&f\end{matrix} \right\}
$$
where $p=(a+b+c+d)/2$. (Since $e,f$ are fixed, one may just as well view this
as a symmetry of the corresponding function $U$.)

More geometrically one may also view this as a symmetry of the set of
tetrahedra in $\IR^3$ (cf. \cite{PzRg, Roberts-6j1}).
First note that if $V_e$ appears in $V_a\otimes V_b$ then the triangle
inequalities 
$$\bigl\vert a-b \bigl\vert \le e \le a+b$$
hold; i.e. there exists a Euclidean triangle with sides of lengths $a,b,e$.
(One also has the parity condition that $a+b+e$ be even.)
Thus if both vectors $v_e$ and $w_f$ are non-zero then there are four 
Euclidean triangles with side lengths $abe,cde,adf,bcf$ respectively.
Now these four triangles may or may not fit together to form the faces of a
Euclidean tetrahedron; the condition that they do is given by 
requiring the determinant of the `Cayley-Menger matrix': 

$$
 \left( \begin {array}{ccccc} 0&{a}^{2}&{e}^{2}&{d}^{2}&1
\\\noalign{\medskip}{a}^{2}&0&{b}^{2}&{f}^{2}&1\\\noalign{\medskip}{e}
^{2}&{b}^{2}&0&{c}^{2}&1\\\noalign{\medskip}{d}^{2}&{f}^{2}&{c}^{2}&0&
1\\\noalign{\medskip}1&1&1&1&0\end {array} \right)
$$
to be positive.
It is simple to check the Regge transformation preserves the set of all 
triangle inequalities (although permuting them in a non-trivial way).
Moreover a computation will show 
that the determinant of the Cayley-Menger matrix is preserved too.

Thus one may view the Regge transformation as an automorphism of the set of
Euclidean tetrahedra, even if we allow real (not necessarily integral) edge
lengths (noting that a tetrahedron with non-zero volume is determined by 
its edge lengths up to isometry, possibly reversing the orientation).

Our first step in deriving the Regge symmetry 
is to note a remarkably similar symmetry of a completely
different object, this time of a certain nonlinear differential equation.
The Painlev\'e VI differential equation (henceforth \PVI) is the following
nonlinear ordinary differential equation for a holomorphic function $y(t)$ with
$t\in\IC\setminus\{0,1\}$:
$$\frac{d^2y}{dt^2}=
\frac{1}{2}\left(\frac{1}{y}+\frac{1}{y-1}+\frac{1}{y-t}\right)
\left(\frac{dy}{dt}\right)^2
-\left(\frac{1}{t}+\frac{1}{t-1}+\frac{1}{y-t}\right)\frac{dy}{dt}  $$
$$
\quad\ +\frac{y(y-1)(y-t)}{t^2(t-1)^2}\left(
\al+\be\frac{t}{y^2} + \gamma\frac{(t-1)}{(y-1)^2}+
\delta\frac{ t(t-1)}{(y-t)^2}\right) $$
where $\al,\be,\ga,\de$  are four complex constants.
\PVI\ is usually thought of as controlling 
the monodromy preserving (isomonodromic) 
deformations of rank $2$ (traceless) Fuchsian systems with $4$ poles on
$\IP^1$ (whose monodromy is a representation of the free group on $3$
generators into $\SL_2(\IC)$).
Okamoto \cite{OkaPVI} proved \PVI\  has a quite nontrivial symmetry:

\begin{thm}\label{thm: ok.sym}
Choose four complex constants $\theta=(\theta_1,\ldots\theta_4)$ and set
$$\al=(\theta_4-1)^2/2, \qquad \be=-\theta_1^2/2, \qquad 
\ga=\theta^2_3/2, \qquad \delta=(1-\theta_2^2)/2.
$$%
If y(t) is a solution of \PVI\ with parameters $\theta$
then, if defined,
$$y+\phi/x$$ solves \PVI\ with parameters
$$\theta'=(\theta_1-\phi,\theta_2-\phi,\theta_3-\phi,\theta_4-\phi)$$
where $\phi=\sum_1^4\theta_i/2$ and
$$
2 x=
\frac {\left( t-1 \right)y'-\theta_1}{y}+
\frac{y'-1-\theta_2}{y-t}-
\frac {t\,y'+\theta_3}{y-1}.
$$
\end{thm}

(Observe the striking similarity with the Regge transformation.)

In the next section 
we will describe exactly how the
Okamoto and Regge symmetries are related (in effect showing precisely how the
complicated Okamoto action on the pair $(y,x)$ relates to the trivial Regge
action on the pair $(e,f)$).

\begin{rmk}
Since the parameters appearing in \PVI\ are now quadratic functions of the
$\theta$'s, if $y$ solves \PVI\ with parameters $\theta$ then $y$ will also
solve \PVI\ for any parameters obtained from $\theta$ by negating any
combination of $\theta_1,\theta_2,\theta_3$ 
and possibly replacing $\theta_4$ by $2-\theta_4$.
Together with the Okamoto transformation these four `trivial' transformations
generate a group isomorphic to the affine Weyl group of type $D_4$ 
(see \cite{OkaPVI}). Further one may add in transformations corresponding to
the $\sym_4$ symmetry group of the affine $D_4$ Dynkin diagram and obtain a
symmetry group isomorphic to the affine Weyl group of type $F_4$. The
confusing fact to note is that one still does
not obtain symmetries corresponding to 
all the tetrahedral $6j$ symmetries, basically because the \PVI\ flows vary
$y,x$ and fix the $\theta$'s.
\end{rmk}

\end{section}

\begin{section}{Regge and Okamoto}

We will relate the Regge action on Euclidean tetrahedra to the 
Okamoto action.

Let $\bfa_1, \bfa_2, \bfa_3\in \IR^3$ be three vectors,
so that ${\bf 0},\bfa_1,\bfa_1+\bfa_2, \bfa_1+\bfa_2+\bfa_3$ 
are the vertices of a tetrahedron.
Denote the other three edge vectors of the tetrahedron by 
$\bfa_4,\bfa_5,\bfa_6$ so that:
$$\bfa_1+\bfa_2+ \bfa_3+\bfa_4=0,$$
$$\bfa_5=\bfa_1+\bfa_2,\qquad\bfa_6=\bfa_2+\bfa_3.$$
Denote the lengths of these six vectors $\bfa_1,\ldots,\bfa_6$ 
by $a,b,c,d,e,f$ respectively.

Now let $\cH$ be the set of traceless $2\times 2$ Hermitian 
matrices. Thus $\cH$ is a real three-dimensional vector space 
and we give it a Euclidean inner product by defining
$$\langle A_1,A_2\rangle := 2\tr(A_1A_2).$$ 

Thus we can view the tetrahedron as living in $\cH$ by choosing 
an isometry $\varphi:\IR^3\cong \cH$ such as
$$
\varphi\left(\begin{matrix} x\\y\\z \end{matrix}\right)=
\frac{1}{2}
\left(\begin{matrix} x & y+iz\\y-iz & -x \end{matrix}\right).
$$
Then we set $$A_j=\varphi(\bfa_j)\in \cH$$ 
for $j=1,\ldots,6$.
Thus the Regge symmetry becomes an action on the set of these $6$ 
Hermitian matrices $A_j$ (clearly determined by its action on the first three
matrices).

Now in the standard isomonodromy interpretation \cite{JM81} 
of the Painlev\'e VI 
equation the Okamoto symmetry becomes a (birational) action on 
the set of Fuchsian systems of the form
\beq\label{eq: lin syst}
\frac{d}{dz}-A, \qquad A:=
\frac{A_1}{z}+\frac{A_2}{z-t}+\frac{A_3}{z-1},\qquad
A_i\in\lsl_2(\IC)
\eeq
where the coefficients $A_1,A_2,A_3$ 
are $2 \times 2$ traceless complex 
matrices and $t\in \IC\setminus \{0,1\}$ fixed.
As above we will 
set $A_4=-(A_1+A_2+A_3)$ which is now the residue of 
$Adz$ at $z=\infty$.
This interpretation comes about by using explicit local coordinates 
on the space of such systems; 
Up to overall conjugation by
$\SL_2(\IC)$, the set of such Fuchsian systems is of complex dimension $6$ 
and local coordinates (near a generic system) are given by 
$$\theta_1,\theta_2,\theta_3,\theta_4,x,y$$ 
where $\theta_i$ is such that $A_i$ has eigenvalues 
$\pm\theta_i/2$, and where $x,y$ are two explicit algebraic functions
of $A$ 
defined for example in \cite{k2p} p.199 (following \cite{JM81}).

Of course one would prefer to view the birational transformation as the
intrinsic object, and its explicit coordinate expression as secondary.
In particular one might hope for a simpler expression than that 
given in terms of
$x,y$ in Theorem \ref{thm: ok.sym}.
One way to do this, which will be useful here,
was observed in \cite{k2p} Lemma 34. 
To describe this we should first 
modify slightly the matrices $A_1,A_2,A_3$: 
Let
$$\wh A_i= A_i+\theta_i/2\qquad i=1,2,3$$
so that $\wh A_i$ has eigenvalues $0,\theta_i$  (i.e. it has rank 
one and trace $\theta_i$).
\begin{lem}[\cite{k2p} Lemma 34]
Each of the five expressions 
$$\tr(\wh A_1\wh A_2),\qquad 
  \tr(\wh A_2\wh A_3),\qquad 
  \tr(\wh A_1\wh A_3),$$
$$ \tr(\wh A_1\wh A_2\wh A_3),\qquad
   \tr(\wh A_3\wh A_2\wh A_1)$$
is preserved by the birational Okamoto transformation of
Theorem \ref{thm: ok.sym}.
\end{lem}

This may be proved by a direct coordinate computation; the geometric 
origin of it is given in \cite{k2p} (see especially Lemma 34, 
Remark 30)
and may be thought of as the
(complexification of the) classical analogue of the ideas we will use 
in the next section.

Now it is straightforward to prove (see \cite{Hit-GASE}) 
that generically 
the first two of these 
expressions (viewed as
functions on the set of Fuchsian systems) 
together with the four $\theta$'s 
make up a system of local coordinates.
Let us write $\la_{12}=\tr(\wh A_1\wh A_2)$ and 
$\la_{23}=  \tr(\wh A_2\wh A_3)$.
Thus, {\em in these coordinates} the Okamoto transformation acts simply as
$$(\theta,\la_{12},\la_{23})\mapsto 
(\theta_1-\phi,\theta_2-\phi,\theta_3-\phi,\theta_4-\phi, 
\la_{12}, \la_{23})
$$
with $\phi=\sum_1^4\theta_i/2$, which looks even more like the 
Regge transformation, and easily yields the main result of this section:

\begin{thm}
Suppose
 $A_1,A_2,A_3$ are $2\times 2$ traceless Hermitian matrices with eigenvalues
$\pm \theta_i/2$ with $\theta_i>0$ 
and let $a,b,c,d,e,f$ be the edge lengths of the
corresponding tetrahedron in $\IR^3$ (under the isometry $\varphi$)
i.e. the lengths of the vectors corresponding to the six matrices
$$A_1,\quad A_2,\quad A_3,\quad
A_1+A_2+ A_3,\quad
A_1+A_2,\quad A_2+A_3$$
respectively.
Then the Okamoto transformation of $(A_1,A_2,A_3)$ corresponds 
to the Regge
transformation of the tetrahedron with edge lengths $a,b,c,d,e,f$.
\end{thm}
\pf
For the first four edge lengths this is easy 
since $(a,b,c,d)=(\theta_1,\theta_2,\theta_3,\theta_4)$ (note that the
triangle inequalities imply $\phi-\theta_i\ge 0$).
We need to also show that the Okamoto transformation preserves the 
edge lengths $e$ and $f$.
But this is now a simple computation: First note
$$e^2=2\tr(A_1+A_2)^2=\theta_1^2+\theta_2^2+4\tr A_1A_2.$$
Then observe
 $4\tr A_1A_2 =4\tr(\wh A_1-\theta_1/2)(\wh A_2-\theta_2/2)
=4\tr\wh A_1\wh A_2-2\theta_1\theta_2$
so that 
$$e^2=(\theta_1-\theta_2)^2+4\tr\wh A_1\wh A_2.$$
The first term on the right here is preserved, as is the second term 
by the lemma above, and so $e$ is preserved since it is positive.
Similarly for $f$.
\epf

\begin{rmk}
Returning briefly to the complex (not-necessarily Hermitian) picture,
the above argument implies that the Okamoto transformation is also
characterised as preserving
$\tr A_5^2$ and $\tr A_6^2$ where $A_5=A_1+A_2$ and 
$A_6=A_2+A_3$ (and similarly it
preserves $\tr (A_1+A_3)^2$).
\end{rmk}

\begin{rmk} (Spherical tetrahedra, cf. \cite{WooTay1}.)
Consider three elements $M_1, M_2, M_3\in \SU(2)$ 
of the $3$-sphere $S^3\cong\SU(2)$, and the spherical tetrahedron with
vertices $I, M_1, M_1M_2, M_1M_2M_3$.
This has edge lengths $l_i$ where $\tr M_i=2\cos(l_i)$ ($i=1,\ldots,6$)
where $M_4=(M_1M_2M_3)^{-1}, M_5=M_1M_2, M_6= M_2M_3$.
One may define a Regge symmetry of the set of such tetrahedra, by acting on
the edge lengths exactly as before.
This action 
may be complexified in the obvious way (allow $M_i\in \SL_2(\IC)$ and
$l_i\in \IC$). 
On the other hand the Okamoto action on the Fuchsian systems 
\eqref{eq: lin syst} induces an action on their monodromy data
(i.e. essentially on the space of $\SL_2(\IC)$ representations of the 
fundamental group of the four-punctured sphere).
The fact to be noted 
is that this action coincides with the above spherical Regge action
(taking $M_i$ to be the monodromy around the $i$th puncture for $i=1,2,3,4$);
in other words the Okamoto action fixes the functions $\tr(M_1M_2)$ and
$\tr(M_2M_3)$ of the monodromy data---this was the main result of  \cite{IIS},
proved differently in Corollary 35 of \cite{k2p}.
\end{rmk}

\end{section}

\begin{section}{Conceptual Regge symmetry}
We will give a conceptual derivation of the fact that the Regge symmetry
preserves the $6j$-symbols (i.e. without using the Racah formula).

There are two basic steps:

$\bullet$ Identify the $\SU(2)$ $6j$ 
symbols with certain $\SU(3)$ $6j$-symbols,

$\bullet$ Use a natural symmetry of these $\SU(3)$ $6j$-symbols. 

These are representation theoretic analogues of the derivation given in
\cite{k2p, srops} of the Okamoto symmetry. (The article \cite{VTL-duke} helped
us to understand this---see also \cite{bafi} where we first realised 
that
\cite{VTL-duke} describes a representation theoretic 
analogue of some things the 
present author had been thinking about.)

First we will set up notation for representations of $H:=\SU(3)$.
Let $W$
be the three-dimensional Hermitian vector space defining $H$.
For any integer $a\ge 0$ write 
$$W_a:=\sym^a(W)$$
for the $a$th symmetric power of $W$; an irreducible representation of $H$.
Similarly for integers $a\ge b\ge 0$ write $W_{(a,b)}$ 
for the irrep of $H$ corresponding to the Young diagram with 3 rows of lengths
$a,b,0$ resp. (cf. \cite{FulHar}). (Thus in particular $W_a=W_{(a,0)}$.)

One would like to define the $\SU(3)$ $6j$-symbols using the same framework 
as described above for $\SU(2)$. 
This is difficult though since $\SU(3)$ is {\em not} 
multiplicity-free and in general one will not obtain a decomposition of the
multiplicity spaces into one-dimensional pieces, but into pieces of higher
dimension. (There are numerous articles discussing this multiplicity
problem, and methods to circumvent it.)
However things are simpler if we take the three initial representations to be
symmetric representations (i.e. of the form $W_a$). Then the Pieri rules
imply 
one will again get the desired one-dimensional decomposition and we may
proceed as before (and this is the only case we will need here).

Thus we choose $3$ symmetric irreps, with labels $a,b,c$ say, and form the 
tensor product
$$W_{abc}:=W_a\otimes W_b\otimes W_c$$
Now, given an arbitrary representation $W_{\la}$ with $\la=(p,q)$ 
one obtains a
multiplicity space as before
$$N_{abc\la}:=\Hom_H(W_\la,W_{abc}).$$

Similarly to before the two expansions
$$W_a\otimes W_b\cong \bigoplus_{(r,s)}
W_{(r,s)}\otimes \Hom_H(W_{(r,s)},W_a\otimes W_b)$$
$$W_b\otimes W_c\cong \bigoplus_{(t,u)}
W_{(t,u)}\otimes \Hom_H(W_{(t,u)},W_b\otimes W_c)$$
yield two decompositions of the multiplicity space:
$$N_{abc\la}\cong \bigoplus_{(r,s)}
 \Hom_H(W_\la,W_{(r,s)}\otimes W_c)\otimes \Hom_H(W_{(r,s)},W_a\otimes W_b)
$$
$$N_{abc\la}\cong \bigoplus_{(t,u)}
 \Hom_H(W_\la,W_a\otimes W_{(t,u)})\otimes \Hom_H(W_{(t,u)},W_b\otimes W_c).
$$
Now the Pieri rules (see \cite{FulHar}) imply that the tensor product of a
symmetric representation with any irrep will be multiplicity free (i.e. each
irrep that appears in the tensor product will appear exactly once).
Thus both of these decompositions of $N_{abc\la}$ will be into one-dimensional
pieces and by choosing real basis vectors of length one $v_{rs}, w_{tu}$ in
the corresponding pieces we can define matrix entries as before:

$$U^{(3)}(a,b,c,\la,(r,s),(t,u)) = \langle v_{rs}, w_{tu} \rangle.$$

(Again a sign-convention is needed to fix the signs of the basis vectors, but
we will not worry about this here---the correct choices will be such that the
following result is true.)

The basic fact that we will use is that these 3-dimensional $6j$ coefficients
are equal to 2-dimensional $6j$ coefficients:

\begin{prop} \label{prop: u=u3}
Choose six integers $a,b,c,d,e,f\ge 0$.
Then [up to sign]
$$U(a,b,c,d,e,f)=U^{(3)}(a,b,c,(p,q),(r,s),(t,u))$$
where 
$$p= \frac{a+b+c+d}{2}, \   r= \frac{a+b+e}{2},\   t= \frac{b+c+f}{2}$$
$$q= p-d,\quad   s= r-e,\quad   u= t-f.$$
\end{prop}

This will be proved below (it is probably quite well-known).
First  we will
describe the desired three-dimensional symmetry and deduce the Regge
symmetry:

\begin{prop} \label{prop: 3dsym}
The three-dimensional $6j$-symbol is symmetric as follows:
$$U^{(3)}(a,\,b,\,c,\,(p,q),\,(r,s),\,(t,u))=$$
$$U^{(3)}(p-a,\, p-b,\,p-c,\,(p, p-q),\,(p-s, p-r),\,(p-u, p-t)).$$
\end{prop}

Note that this does indeed project back to give the Regge symmetry\footnote{
Inverting
the equations appearing in Proposition \ref{prop: u=u3} yields
$d=p-q,e=r-s,f=t-u$ so that Proposition \ref{prop: u=u3} implies
$U^{(3)}(p-a,\,p-b,\,p-c,\,(p,\, p-q),\,(p-s,\, p-r),\,(p-u,\, p-t))$
$=U(p-a,\,p-b,\,p-c,\,q,\,r-s,\,t-u)$
and so, since $q=p-d,\, r-s=e,\, t-u=f$ 
we do obtain the Regge symmetry
$U(a,\,b,\,c,\,d,\,e,\,f)=U(p-a,\,p-b,\,p-c,\,p-d,\,e,\,f)$ as desired.}.
Thus our task is reduced to justifying the above two propositions.
At first sight it may appear that little progress has been made, replacing the
Regge symmetry by the symmetry of Proposition \ref{prop: 3dsym}. But as we
shall see, this symmetry arises simply by pairing two dual representations of
$\SU(3)$. [It is not however simply a matter of dualising all the
  representations in sight, since 
the dual of $W_a$ is not a
  symmetric representation, for $a>0$.]

\pfms (of Proposition \ref{prop: u=u3}).
The first step is to identify the corresponding weight spaces $M_{abcd}$ and
$N_{abc\la}$.
This can be done easily using ``Howe duality for $\GL_k$-$\GL_n$'', as follows
(cf. \cite{Zhelo,VTL-duke})\footnote{To proceed explicitly (using similar
  ideas) see \cite{KlinkT, GliskeKT}.}.
Choose two positive integers $k,n$ and let $V,W$ be complex vector spaces of
dimensions $k,n$ respectively.
Then their tensor product $V\otimes W$ is a representation of 
$\GL(V)\times \GL(W)$
and so  its $d$th symmetric power $\sym^d(V\otimes W)$ is also a 
$\GL(V)\times \GL(W)$-module.
This decomposes as a direct sum of irreducible 
$\GL(V)\times \GL(W)$-modules in the following way (see \cite{FulHar}
Exercise 6.11; apparently this goes back to Cauchy):
$$\sym^d(V\otimes W)\cong 
\bigoplus_{
\begin{smallmatrix}
\#\la\le k,n\\ \vert\la\vert=d
\end{smallmatrix}}
V_\la\otimes W_\la$$
where the sum is over all Young diagrams $\lambda$ 
with $d$ boxes and having no more than
$k$ or $n$ rows, and $V_\la$ (resp. $W_\la$) is the irreducible 
$\GL(V)$-module (resp. $\GL(W)$-module) corresponding to $\lambda$.
Choosing bases of $V$ and $W$ allows us to be more explicit.
In particular it identifies $\GL(V)\cong\GL_k(\IC)$ and so picks out a maximal
torus (the diagonal subgroup),
as well as a Borel subgroup (the upper triangular subgroup) and so allows us
to speak of weights and highest weight vectors of $\GL(V)$ modules (similarly
for $\GL(W)$).
Also we can now view $V\otimes W$ as the space of linear functions $\psi$
on the set $M_{k\times n}$ of $k\times n$ matrices $X=(x_{ij})$ 
with the following action of $\GL_k(\IC)\times \GL_n(\IC)$: 
$$(g_k,g_n)(\psi)(X)=\psi(g_k^TXg_n).$$
(To avoid confusion below when $k=n=3$, we will refer to this
$\GL_k(\IC)$-action, respectively $\GL_n(\IC)$-action,
as the action on the left, resp. right.)
Then $\sym^d(V\otimes W)=\sym^d(M_{k\times n}^*)$ 
is the set of such functions which are homogeneous polynomials of
degree $d$, with the same action.
Thus as $\GL_k(\IC)\times \GL_n(\IC)$-modules
\beq \label{eq: exp decomp.}
\sym^d(M_{k\times n}^*)\cong 
\bigoplus_{
\begin{smallmatrix}
\#\la\le k,n\\ \vert\la\vert=d
\end{smallmatrix}}
V^{(k)}_\la\otimes V^{(n)}_\la
\eeq
where $V_\la^{(r)}$ is the $\GL_r(\IC)$ irrep. with Young diagram $\la$.

Now, as a $\GL_k$ module $\sym^\bullet(M_{k\times n}^*)$ can be viewed as the
tensor product of the functions on each of the columns of $X$, i.e.
$\sym^\bullet(M_{k\times n}^*)\cong \bigotimes_1^n\Sym^\bullet\IC^k$. 
Moreover if
we choose $n$ positive integers $\mu=(\mu_1,\ldots,\mu_n)$ then we can
consider the subspace:
$$S^\mu\IC^k:=\sym^{\mu_1}\IC^k\otimes\cdots\otimes\sym^{\mu_n}\IC^k\subset 
\sym^\bullet(M_{k\times n}^*)$$
of functions which are homogeneous of degree $\mu_i$ in the $i$th column.
This is equivalent to saying they are the vectors of weight $\mu$ for the 
$\GL_n$ action on $\sym^\bullet(M_{k\times n}^*)$.

Now what we are really interested in are the $\GL_k$ 
multiplicity spaces of the form:
$$M_\la^\mu:=\Hom_{\GL_k}(V_\la^{(k)},S^\mu\IC^k)$$
whose dimension is the multiplicity of $V_\la^{(k)}$ 
in this $n$-fold product of 
symmetric representations.
This
 multiplicity space
may be realised explicitly as the subspace of $S^\mu\IC^k$ of vectors of
highest weight $\la$ for the $\GL_k$ action 
(since each copy of $V_\la^{(k)}$ in $S^\mu\IC^k$ has a unique highest
weight vector, up to scale).
On the other hand by the decomposition 
\eqref{eq: exp decomp.} the subspace of $\sym^\bullet(M_{k\times n}^*)$
of vectors of
highest weight $\la$ for the $\GL_k$ action
is a single copy of $V_\la^{(n)}$. Intersecting this with the subspace
$S^\mu\IC^k$ (i.e. the vectors with $\GL_n(\IC)$ weight $\mu$) 
yields the basic result we need (cf. \cite{VTL-duke} Lemma 3.4):
\begin{lem}\label{lem: wt-mult im}
The above discussion gives an isomorphism  
$V_\la^{(n)}[\mu]\cong M^\mu_\la$
between the weight $\mu$ subspace of the $\GL_n(\IC)$ 
representation $V_\la^{(n)}$
 and the $\GL_k(\IC)$ multiplicity space  $M^\mu_\la$.
\end{lem}

In the $6j$-symbol situation, 
we are interested in three-fold tensor products of symmetric
representations, so $n=3$, and our Young diagrams always have at most two
non-zero rows, {\em so we may take any $k\ge 2$}.
For $k=2$, Lemma \ref{lem: wt-mult im} implies
$$M_{abcd}\cong V^{(3)}_\la[\mu]$$
where $\mu=(a,b,c)$ and where $\la$ has at most $2$ rows and 
is such that $V^{(2)}_\la\cong V_d$ as $\SU(2)$ representations  
and the number of boxes in $\la$ equals $a+b+c$. 
This implies $\la=(p,q,0)$ where $p+q=a+b+c$ and $p-q=d$, i.e.
$p=(a+b+c+d)/2$ and $q=p-d$ as in the statement of Proposition 
\ref{prop: u=u3}.
In turn, for $k=3$, Lemma \ref{lem: wt-mult im} implies
$$V^{(3)}_\la[\mu] \cong N_{abc\la}$$
and so combining the two gives the desired isomorphism of multiplicity spaces.
(Explicitly if we view the $2\times 3$ matrices as the first two rows
of the $3\times 3$ matrices, then these 
multiplicity spaces are actually equal as spaces of 
polynomial functions on $3\times 3$ matrices since in both
cases they do not depend on the variables in the third row, as $\lambda$
has at most two non-zero rows.)

\begin{rmk}
For the reader familiar with Gelfand-Tsetlin tableaux, we should mention 
that the weight space $V^{(3)}_\la[\mu]$ 
(and thus the common multiplicity
space)
admits a basis parameterised by tableaux of the form
$$
\left(\begin{matrix}
p &&   q  && 0 \\
& \al && \be & \\
&&    \ga &&
\end{matrix}\right)$$
where $\al,\be,\ga$ are integers satisfying interlacing inequalities:
 $p\ge\al\ge q\ge\be\ge 0,\ \al\ge\ga\ge\be$
and should be such that the tableau has `weight' $\mu=(a,b,c)$---the 
weight of a tableau is the differences of the row sums, i.e. we require
$\ga=a,\ \al+\be-\ga=b,\  p+q -(\al+\be)=c.$
This gives a simple way to compute the dimension of the multiplicity spaces
(i.e. count the tableaux) although we will not need to 
use this Gelfand-Tsetlin basis
(in general it does not coincide with any of the  coupling bases).
\end{rmk}

To complete the proof of Proposition \ref{prop: u=u3} 
we need to see the corresponding coupling bases match up
under the above isomorphism of multiplicity spaces.
(One may use the Bargmann-Segal-Fock Hermitian form on 
$\sym^\bullet(M_{k\times n}^*)$ 
and so  see the Hermitian forms coincide.)
One way to do this is to first observe that the $1$-$2$
coupling spaces are the eigenspaces of the $\GL_k(\IC)$ 
quadratic Casimir operator $\cC_{k}$
acting on the first two tensor factors of 
$S^\mu\IC^k=\sym^a\IC^k\otimes \sym^b\IC^k\otimes \sym^b\IC^k$, where
$\mu=(a,b,c)$ and $k=2$ or $3$.
This holds since in general $\cC_{k}$ acts 
(see e.g. \cite{Zhelo} p.161)
by multiplication by the scalar 
\beq\label{eq: cas action}
\sum m_i^2+\sum_{i<j}m_i-m_j\eeq on the $\GL_k$ 
irrep. with Young diagram $(m_1,\ldots,m_k)$.
Denote this operation by $\cC^{12}_k$ and note
it preserves the common multiplicity space 
$$M_{abcd}=N_{abc\la}\subset S^\mu\IC^2\subset S^\mu\IC^3$$
for $k=2,3$.
Then we just observe, using \eqref{eq: cas action}, 
that on the multiplicity space
the Casimirs differ by a scalar:
$\cC^{12}_3=\cC_2^{12}+ (a+b)$
and so have the same eigenspaces (and the eigenspace
labels match up as stated, that is: $r-s=e,r+s=a+b$).
(Similarly for the $2$-$3$ coupling.)

\epfms

Next we will describe the three-dimensional symmetry which lifts the
Regge symmetry: 
If $V^{(3)}_\la$ 
is the irrep of $\GL_3(\IC)$ with Young diagram $\la=(p,q,0)$ 
then there is a pairing
$$V^{(3)}_\la\otimes (V^{(3)}_\la)^\vee\to \IC$$
where $(V^{(3)}_\la)^\vee$ is the dual representation. Tensoring this by the $p$th
power $D^p$ of the determinant representation yields a pairing
$V^{(3)}_\la\otimes V^{(3)}_{\la'}\to D^p$, where $\la'=(p,p-q,0)$.
This pairs the $\mu=(a,b,c)$ 
weight space of $V^{(3)}_\la$ with the $\mu'$ weight space 
of $V^{(3)}_{\la'}$ where $\mu'=(p-a,p-b,p-c)$ and so yields a perfect pairing:
\beq \label{eq: pairing on wt spaces}
V^{(3)}_\la[\mu]\otimes V^{(3)}_{\la'}[\mu']\to \IC.\eeq
Now to prove Proposition \ref{prop: 3dsym} one just 
needs to check that the
corresponding coupling bases are dual with respect to this pairing.

\begin{rmk}
On the level of Gelfand-Tsetlin tableaux the above duality corresponds to
negating each tableau element then adding $p$ to each element and finally 
flipping the tableau about its vertical axis. (Observe that the
tableau's weight has transformed as stated.)
\end{rmk}

\pfms (of Proposition \ref{prop: 3dsym}).
We will first 
describe the above pairing in a different way which will be more convenient.
Write $G=\GL_3(\IC)$ and let $S^a=\sym^a\IC^3$. 
By the Pieri rules there is a unique $G$-equivariant map
$$S^a\otimes S^{p-a}\to S^p$$
since $S^p$ appears precisely once in this tensor product\footnote{
More precisely there is a unique subspace of $S^aS^{p-a}$ which is 
$G$-equivariantly 
isomorphic to $S^p$,
and thus a (unique) orthogonal projection onto this subspace. To lighten the
notation we will call this subspace  $S^p$ here (and similarly below).}.
(Similarly replacing $a$ by $b$ or $c$.)
Putting these together there is a $G$-equivariant map
$$S^aS^bS^c\otimes S^{p-a}S^{p-b}S^{p-c}\to S^pS^pS^p$$
pairing the corresponding factors (and omitting to write several $\otimes$
symbols). Again by the Pieri rules there is a unique (projection) map
$ S^pS^pS^p\to D^p$ to the $p$th power of the determinant representation.
Composing with the above map we get a (degenerate) pairing:
$$\nu:S^aS^bS^c\otimes S^{p-a}S^{p-b}S^{p-c}\to D^p.$$
In terms of pairs of 
polynomials on $M_{3\times 3}$ (viewing $S^\al S^\be S^\ga$ as polynomials
homogeneous of degrees $\al,\be,\ga$ in the columns $1,2,3$ resp. as above)
this bilinear form $\nu$ 
amounts to multiplication followed by orthogonal
projection onto $D^p$ (which 
is just the one dimensional subspace spanned $p$th power of the polynomial 
$\det:M_{3\times 3}\to \IC$).

Now we wish to relate $\nu$ to the natural pairing
$V^{(3)}_\la[\mu]\otimes V^{(3)}_{\la'}[\mu']\to \IC$
of \eqref{eq: pairing on wt spaces}.
For this we view $V^{(3)}_\la$ as a space of polynomials on the $3\times 3$
matrices as in Lemma \ref{lem: wt-mult im} (as the polynomials 
with highest weight 
$\la$ for the $\GL_3$-action on the left),
but we view $V^{(3)}_{\la'}$
differently as the space of polynomials with {\em lowest} weight 
$(0,p-q,p)$ for the $\GL_3$-action on the left.
Then multiplication of functions followed by orthogonal
projection onto $D^p$ gives a pairing 
$V^{(3)}_\la\otimes V^{(3)}_{\la'}\to D^p$
which one may check is nonzero directly (observing that
a highest weight vector of 
$V^{(3)}_\la$ pairs non-degenerately
with a lowest weight vector of $V^{(3)}_{\la'}$).
Thus by Schur's lemma this pairing coincides with 
the natural one up to scale.
Restricting to the $\mu$ and $\mu'$ weight spaces respectively and using Lemma
 \ref{lem: wt-mult im}
thus yields a non-degenerate pairing
\beq\label{eq: wt space pairing}
N_{\mu\la}\otimes N_{\mu'\la'}\to \IC
\eeq
which is a restriction of $\nu$.
(To identify $N_{\mu'\la'}\cong V^{(3)}_{\la'}[\mu']$ we use the analogue of 
Lemma  \ref{lem: wt-mult im} with ``highest weight'' replaced by
``lowest weight'' throughout the proof.)  
Now we need to show that with respect to the pairing 
\eqref{eq: wt space pairing} the 
$1$-$2$ coupling bases on each side, are dual (and
   similarly for the $2$-$3$ bases). 
For this it is sufficient to prove the 
non-corresponding coupling basis vectors are orthogonal 
(since we know the pairing is non-degenerate this forces the corresponding
coupling basis vectors to pair up).
We will show this for the $1$-$2$ coupling (the other coupling being
analogous).

Write $W_{(x,y)}$ for the irrep. of $G$ with Young diagram
$(x,y,0)$. By the Pieri rules there is a unique
map $S^pS^p\to W_{(p,p)}$ and so we have a $G$-equivariant map
\beq\label{eq: smallpair}
S^aS^bS^{p-a}S^{p-b}\to S^pS^p\to W_{(p,p)}.
\eeq
This enables us to factor $\nu$ as follows:
$$S^aS^bS^cS^{p-a}S^{p-b}S^{p-c}\to W_{(p,p)}S^cS^{p-c}\to 
W_{(p,p)}S^p\to D^p.$$
To see how the coupling subspaces pair up, 
first expand both $S^aS^b$ and $S^{p-a}S^{p-b}$ into sums of irreps 
so $S^aS^bS^{p-a}S^{p-b}$ becomes
a sum of tensor products of the form $W_{(r,s)}\otimes W_{(x,y)}$, where
$W_{(r,s)}\subset S^aS^b$ etc., and \eqref{eq: smallpair} maps these to 
 $W_{(p,p)}$. However using the Littlewood--Richardson rule it is easy to see
that there is a nonzero map 
$W_{(r,s)}\otimes W_{(x,y)}\to W_{(p,p)}$ if and only if
$x=p-s,y=p-r$ (and if so it is unique up to scale). This gives the
stated correspondence between the $1$-$2$ coupling bases (in the 
fifth slot of $U^{(3)}$).

\epfms

\begin{rmk}
Presumably the argument above extends to the Regge symmetry of the
q-deformation of the  $6j$-coefficients \cite{KirResh6j}.
\end{rmk}

\end{section}

\renewcommand{\baselinestretch}{0.9}              %
\normalsize
\bibliographystyle{amsplain}    \label{biby}
\bibliography{../thesis/syr}    
\end{document}